\documentclass[12pt,a4paper]{article}
\usepackage{latexsym,amssymb,amsfonts,amsbsy,amsmath,bm}
\usepackage{graphics}
\usepackage[utf8]{inputenc}
\usepackage[margin=2.5cm]{geometry}
\usepackage{colortbl}
\usepackage{amsthm}
\usepackage{hyperref}

\newcommand{\BigO}{\mathcal{O}}
\title{Computations for  the first Lyapunov coefficient.}
\author{
	Marino Badiale\thanks{Department of Mathematics, University of Turin, Email: \href{mailto:marino.badiale@unito.it}{marino.badiale@unito.it}} 
	\and
	Isabella Cravero\thanks{Department of Mathematics, University of Turin, Email: \href{mailto:isabella.cravero@unito.it}{isabella.cravero@unito.it}. Corresponding author.}  }
\date{}
\begin{document}
	\maketitle

	These notes are a supplementary file to the paper~\cite{BaCra_5}, 
	where we present in full detail the computations developed in Section~4.2. 
	The purpose of that section is to compute the first Lyapunov coefficient 
	\(a(\mu(\varepsilon),\varepsilon)\), given according to 
	Kuznetsov~\cite[Ch.~3, Sec.~5]{KuznetsovFourth} by
	\begin{equation} \label{eq:Lyapcoeff}
		a(\mu(\varepsilon),\varepsilon ) 
		= \frac{1}{2\omega_0} \operatorname{Re} \left[
		\langle p, C(q, q, \bar{q}) \rangle 
		- 2 \langle p, B(q, A^{-1} B(q, \bar{q})) \rangle 
		+ \langle p, B(\bar{q}, (2i\omega_0 I_n - A)^{-1} B(q, q)) \rangle
		\right],
	\end{equation}
	focusing on its asymptotic developments with respect to~\(\varepsilon\). 
	
	Here \(A\) denotes the Jacobian matrix of the vector field \(F\) evaluated at the equilibrium point 
	\(P_0 = P_0(\mu(\varepsilon), \varepsilon)\). 
	Along the curve \(\mu = \mu(\varepsilon)\), the matrix \(A\) has a pair of purely imaginary conjugate eigenvalues 
	\(\pm i \omega_0\). 
	This curve admits a Taylor expansion for \(\varepsilon \to 0\), starting from the value 
	\(\mu(0) = \mu_0\), which corresponds to the unperturbed case \(\varepsilon = 0\). 
	
	The vector \(q \in \mathbb{C}^4\) is a right eigenvector of \(A\) associated with \(i \omega_0\),
	and \(p \in \mathbb{C}^4\) is the corresponding left eigenvector, normalized so that 
	\(\langle p, q \rangle = 1\).
	
	The multilinear operators \(B(x, y)\) and \(C(x, y, z)\) are defined in terms of the derivatives of the vector field 
	\(F(\xi_1, \xi_2, \xi_3, \xi_4) = (F_1, F_2, F_3, F_4)\) evaluated at \(P_0\):
	\begin{equation}\label{eq:B}
		B_i(x, y) 
		= \sum_{j,k=1}^{4} 
		\left. 
		\frac{\partial^2 F_i(\xi)}{\partial \xi_j \partial \xi_k} 
		\right|_{\xi = P_0} 
		x_j y_k,
		\qquad i = 1, \dots, 4,
	\end{equation}
	and
	\begin{equation}\label{eq:C}
		C_i(x, y, z) 
		= \sum_{j,k,l=1}^{4} 
		\left. 
		\frac{\partial^3 F_i(\xi)}{\partial \xi_j \partial \xi_k \partial \xi_l} 
		\right|_{\xi = P_0} 
		x_j y_k z_l,
		\qquad i = 1, \dots, 4.
	\end{equation}
	Here \(B\) and \(C\) denote the symmetric bilinear and trilinear forms corresponding 
	to the second and third derivatives of the vector field \(F\), respectively.
	Their explicit expressions will be computed in the following sections.
	
	\vskip 0.3cm
	\noindent
	The exposition follows the same scheme as in Section~3.2 of the paper~\cite{BaCra_5} 
	and the algorithm described in~\cite[Ch.~3, Sec.~5]{KuznetsovFourth}.

	\vskip 0.3cm
	\noindent
	We start from the model studied in the main paper~\cite{BaCra_5}, given by
	\begin{equation} \label{eq:modello}
		\left\{
		\begin{aligned}
			x^\prime &= -2x + 3w, \\
			y^\prime &= 50 y - \tfrac{1}{2} y^2 - \delta x y, \\
			z^\prime &= - \varepsilon x z + k \varepsilon \frac{xw}{\varepsilon x w +1} z,  \\
			w^\prime &= \varepsilon x z - \varepsilon x w - 5 w + \delta x y.
		\end{aligned}
		\right.
	\end{equation}
	
	\noindent
	Setting \( \mu = \tfrac{5}{3 \delta} \), we obtain the equilibrium point 
	\( P_0(\mu, \varepsilon) = (x_0, y_0, z_0, w_0)\), given by
	\[
	P_0(\mu, \varepsilon) =\left(
	\frac{15 \mu (50 - \mu)}{25 + 3 \varepsilon \mu^2}, \;
	\frac{50 \mu (6 \varepsilon \mu + 1)}{25 + 3 \varepsilon \mu^2}, \;
	0, \;
	\frac{10 \mu (50 - \mu)}{25 + 3 \varepsilon \mu^2}
	\right),
	\]
	whose Jacobian matrix at the equilibrium point \(P_0(\mu, \varepsilon)\) is given by
	\begin{equation} \label{eq:jacobiano_punto_v2}
		J(P_0(\mu, \varepsilon)) = \begin{pmatrix} 
			-2 & 0 & 0 & 3 \\[6pt]
			\displaystyle - \tfrac{10}{3} + \varepsilon d_1 & - \mu + \varepsilon d_2 & 0 & 0 \\[6pt]
			0 & 0 & \eta_0 & 0 \\[6pt]
			\displaystyle \tfrac{10}{3}  & 50 - \mu + \varepsilon d_2 & \varepsilon x_0 & -5 - \varepsilon x_0
		\end{pmatrix},
	\end{equation}
	where
	\[
	\begin{aligned}
		d_1 & = d_1(\mu, \varepsilon) =     10  \frac{\mu - 50}{25 + 3 \varepsilon \mu^2}  \mu,  \\
		d_2 & = d_2(\mu, \varepsilon) =   3  \frac{\mu - 50}{25 + 3 \varepsilon \mu^2}  \mu^2, \\
		\displaystyle \eta_0 & = - \varepsilon x_0  + k \frac{\varepsilon x_0 w_0}{\varepsilon x_0 w_0 +1},
	\end{aligned}
	\]
	 with \(\displaystyle k < \frac{3}{2 x_0}\) (recalling that \(w_0 = \frac{2}{3} x_0\)), so that \(\eta_0 < 0\).
	Here the functions \(d_i(\mu,\varepsilon), i=1,2\), are \(C^\infty\) in a neighborhood of \(\mu = \mu_0\) (the critical value introduced in the main paper, 
	to be recalled later in this note) and \(\varepsilon = 0\).
	
	\noindent
	It is immediate that \(\eta_0\) is an eigenvalue of \(J(P_0(\mu, \varepsilon))\). The other three eigenvalues coincide with those of the reduced Jacobian
	\begin{equation} \label{eq:jacobiano_ridotto}
		J_1(P_0(\mu, \varepsilon)) = \begin{pmatrix} 
			-2 & 0 & 3 \\[6pt]
			\displaystyle - \tfrac{10}{3} + \varepsilon d_1 & - \mu + \varepsilon d_2 & 0 \\[6pt]
			\displaystyle \tfrac{10}{3}  & 50 - \mu + \varepsilon d_2 & -5 - \varepsilon x_0
		\end{pmatrix}.
	\end{equation}
	To evaluate \(a(\mu(\varepsilon),\varepsilon)\) we need asymptotic expansions of the entries of \eqref{eq:jacobiano_ridotto}. 
	First, we compute
	\[
	\frac{1}{25 + 3 \varepsilon \mu^2}= \frac{1}{25} - \frac{3}{625} \mu^2 \varepsilon + \BigO(\varepsilon^2),
	\]
	so that
	\[
	\begin{aligned}
		d_1 &=  \frac{2}{5} \mu (\mu-50) + \BigO(\varepsilon), \\
		d_2 &  = \frac{3}{25} \mu^2 (\mu-50) + \BigO(\varepsilon).
	\end{aligned}
	\]
	The reduced Jacobian \eqref{eq:jacobiano_ridotto} is now given by 
	\begin{equation} \label{eq:jacobiano_ridottosviluppo}
		J_1(P_0(\mu, \varepsilon)) = \scalebox{0.9}{$ \begin{pmatrix} 
			-2 & 0 & 3 \\[6pt]
			\displaystyle - \frac{10}{3} + \frac{2}{5} \mu (\mu -50) \varepsilon + \BigO(\varepsilon^2) 
			& \displaystyle -\mu + \frac{3}{25} \mu^2 ( \mu -50) \varepsilon + \BigO(\varepsilon^2) 
			& 0 \\[6pt]
			\displaystyle \frac{10}{3} 
			& 50 - \mu + \tfrac{3}{25} \mu^2 ( \mu -50) \varepsilon + \BigO(\varepsilon^2) 
			& -5 + \frac{3}{5} \mu (\mu - 50) \varepsilon + \BigO(\varepsilon^2)
		\end{pmatrix}.$}
	\end{equation}
	Furthermore, the eigenvalue \(\eta_0\) admits the expansion
	\[
	\begin{aligned}
		\eta_0 
		& = \frac{3}{5} \mu (50 - \mu) \left( -1 + \frac{2}{5} k \mu (50 - \mu) \right) \varepsilon + \BigO(\varepsilon^2) \\
		& = x_T \left(-1 + \frac{2}{3} k x_T \right) \varepsilon + \BigO(\varepsilon^2) 
		= \alpha_T \varepsilon + \BigO(\varepsilon^2),
	\end{aligned}
	\]
	where \( x_T = \frac{3}{5} \mu (50 - \mu) \) is the value obtained for the unperturbed  system in the main paper and \(\alpha_T = \left( -1 + \frac{2}{3} k x_T \right) x_T\).

\vskip 0.2cm
\noindent	
The characteristic polynomial of \(J_1(P_0(\mu, \varepsilon))\), computed using the Taylor expansions, becomes
\begin{equation} \label{eq:polcar2}
	\begin{aligned}
		p_{\mu}(x) &= x^3
		 +  \left[7 + \mu - \frac{3}{25} \mu \left( \mu^2 - 45 \mu - 250 \right)  \varepsilon + \BigO(\varepsilon^2) \right] x^2 + \\
		& +  \left[ 7 \mu - \frac{6}{25} \mu \left( 6 \mu^2 - 295 \mu - 250 \right) \varepsilon + \BigO(\varepsilon^2) \right] x + \\
		& + 500 - 10 \mu + \frac{6}{5} \mu (\mu - 50)^2 \varepsilon + \BigO(\varepsilon^2).
	\end{aligned}
\end{equation}

We define
\[
\begin{aligned}
	b_1(\mu) &= -\frac{3}{25} \mu \left( \mu^2 - 45 \mu - 250 \right), \\
	b_2(\mu) &= -\frac{6}{25} \mu \left( 6 \mu^2 - 295 \mu - 250 \right), \\
	b_3(\mu) &= \frac{6}{5} \mu (\mu - 50)^2,
\end{aligned}
\]
so that  
\[
p_{\mu}(x) = x^3 + [7 + \mu + b_1(\mu) \varepsilon + \BigO(\varepsilon^2)] x^2 + [7 \mu + b_2(\mu) \varepsilon + \BigO(\varepsilon^2)] x + 500 - 10 \mu + b_3(\mu) \varepsilon + \BigO(\varepsilon^2).
\]
We recall that the Jacobian matrix has two purely imaginary conjugate eigenvalues 
\(\pm i \omega_0\) provided that its characteristic polynomial 
\(p_{\mu}(x) = x^3 + a_1 x^2 + a_2 x + a_3\) satisfies \(a_1 a_2 - a_3 = 0, \; a_2 >0\), 
with \(\omega_0 = \sqrt{a_2}\).

\noindent
If we consider \(p_{\mu}(x)\) as in \eqref{eq:polcar2}, since \(\mu>0\), the leading term is positive, and hence the quantity \(a_2\) 
is positive for sufficiently small \(\varepsilon\).
 Thus, we set
\[
F(\mu, \varepsilon) = a_1 a_2 -a_3 = 7 \mu^2 + 59 \mu - 500 + [7 \mu b_1(\mu) + (7 + \mu) b_2(\mu) - b_3(\mu)] \varepsilon + \BigO(\varepsilon^2),
\]
and we know that \(F(\mu,0)=7 \mu^2 + 59 \mu - 500=0\) if \(\mu=\mu_0= \frac {-59 + \sqrt{17481}} {14}\). 
By the Implicit Function Theorem, we can therefore determine a curve 
\(\mu(\varepsilon)\) along which the characteristic polynomial admits a pair of complex conjugate roots.

%
%

We already know that \(\mu(0)=\mu_0\) and the derivative of \(\mu(\varepsilon)\) evaluated at \(\varepsilon = 0\) is
\[
\mu_1 = \mu'(0) = \left. -\frac{\partial F / \partial \varepsilon}{\partial F / \partial \mu} \right|_{\substack{\varepsilon = 0 \\ \mu \approx \mu_0}} 
= -\frac{7 \mu_0 b_1(\mu_0) + (7 + \mu_0) b_2(\mu_0) - b_3(\mu_0)}{14 \mu_0 + 59}.
\]

We can therefore conclude that the characteristic polynomial admits a pair of purely imaginary conjugate roots as \(\varepsilon\) varies, for values of \(\mu\) given by
\begin{equation} \label{eq:mu}
	\mu(\varepsilon) = \mu_{\varepsilon} = \mu_0 + \mu_1 \varepsilon + \BigO(\varepsilon^2).
\end{equation}
Therefore, along \(\mu= \mu(\varepsilon)\), we have
\[
 a_2= 7 \mu_{\varepsilon} + b_2(\mu_{\varepsilon}) \varepsilon + \BigO(\varepsilon^2) 
= 7 \mu_0 + (7 \mu_1 + b_2(\mu_0)) \varepsilon + \BigO(\varepsilon^2) 
= 7 \mu_0 + \mu_2 \varepsilon + \BigO(\varepsilon^2)
\]
since \(b_2(\mu_{\varepsilon}) = b_2(\mu_0) + \BigO(\varepsilon)\), with \(\mu_2 = 7 \mu_1 + b_2(\mu_0) \), hence
\[
\omega_0 = \sqrt{7 \mu_0 + \mu_2 \varepsilon + \BigO(\varepsilon^2)} 
= \sqrt{7 \mu_0} + \frac{\mu_2}{2 \sqrt{7 \mu_0}} \varepsilon + \BigO(\varepsilon^2) 
= \sqrt{7 \mu_0} + \mu_3 \varepsilon + \BigO(\varepsilon^2)
\]
with \(\displaystyle \mu_3 = \frac{\mu_2}{2 \sqrt{7 \mu_0}}\).
 So the pair of complex eigenvalues is given by
\begin{equation} \label{eq:rho}
	\rho_{1,2} = \pm i \omega_0, \qquad \omega_0=\sqrt{7 \mu_0} + \mu_3 \varepsilon + \BigO(\varepsilon^2).
\end{equation}
From now on, we will always consider \(\mu= \mu(\varepsilon)\) given by \eqref{eq:mu}, with \(P_0=P_0(\mu(\varepsilon),\varepsilon)\) and therefore we will be in the case where both the  Jacobian \(J(P_0)\)
and the reduced Jacobian \(J_1(P_0)\)
have a pair of purely imaginary complex conjugate eigenvalues given by \eqref{eq:rho}.
%
Then,  we rewrite \(A=J(P_0)\) in \eqref{eq:jacobiano_punto_v2} as
\begin{equation} \label{eq:matriceA}
	A=  \begin{pmatrix} 
		-2 & 0 & 0 & 3 \\ \\
		- \frac{10}{3} + \beta_1 \varepsilon + \BigO(\varepsilon^2) 
		& -\mu_0 + \beta_2 \varepsilon + \BigO(\varepsilon^2) & 0 & 0 \\ \\
		0 & 0 & \alpha_T \varepsilon + \BigO(\varepsilon^2) & 0 \\ \\
		\frac{10}{3}  
		& 50 - \mu_0 + \beta_2 \varepsilon + \BigO(\varepsilon^2) 
		& \beta_3 \varepsilon + \BigO(\varepsilon^2) 
		& -5 - \beta_3 \varepsilon + \BigO(\varepsilon^2)
	\end{pmatrix}
\end{equation}
with
\[
\begin{aligned}
	\beta_1 &=  \frac{2}{5} \mu_0 ( \mu_0 -50), \\
	\beta_2 &= -\mu_1 + \frac{3}{25} \mu_0^2 ( \mu_0 -50), \\
	\beta_3 &= \frac{3}{5} \mu_0 (50 - \mu_0).
\end{aligned}
\]

\noindent
We now search for the complex eigenvectors associated with the eigenvalue \(i \omega_0\). We aim to find \(\varphi = (\varphi_1, \varphi_2, \varphi_3, \varphi_4) \in \mathbb{C}^4\) such that
\begin{equation} 
\scalebox{0.8}{$	\begin{pmatrix} 
		-2-i \omega_0 & 0 & 0 & 3 \\ \\
		- \frac{10}{3} + \beta_1 \varepsilon + \BigO(\varepsilon^2) 
		& -\mu_0 + \beta_2 \varepsilon - i \omega_0+ \BigO(\varepsilon^2) & 0 & 0 \\ \\
		0 & 0 & \alpha_T \varepsilon - i \omega_0+ \BigO(\varepsilon^2) & 0 \\ \\
		\frac{10}{3}  
		& 50 - \mu_0 + \beta_2 \varepsilon + \BigO(\varepsilon^2) 
		& \beta_3 \varepsilon + \BigO(\varepsilon^2) 
		& -5 - \beta_3 \varepsilon - i \omega_0 + \BigO(\varepsilon^2)
	\end{pmatrix} $} \begin{pmatrix}  \varphi_1 \\  \varphi_2 \\  \varphi_3 \\  \varphi_4 \end{pmatrix}= \begin{pmatrix} 0 \\  0 \\  0 \\ 0 \end{pmatrix}
\end{equation}
Solving the system, we find
\begin{equation}
	\left\{
	\begin{aligned}
		\varphi_2 & =  \frac {- 10 + 3 \beta_1 \varepsilon } {3(\mu_0 - \beta_2 \varepsilon + i \omega_0)} \varphi_1 + \BigO(\varepsilon^2), \\
		\varphi_3 & = 0, \\
		\varphi_4 & =  \frac {2 + i \omega_0} {3} \varphi_1.
	\end{aligned}
	\right.
\end{equation}
with \(\varphi_1 \in \mathbb{C}\). 
Considering \(\varphi_1 = 3\left[\mu_0 + i \sqrt{7 \mu_0} + (-\beta_2 + i \mu_3) \varepsilon + \BigO(\varepsilon^2)\right]\), we obtain:
\[
\begin{aligned}
	{q}_1 &= 3[\mu_0 + i \sqrt{7 \mu_0} + (-\beta_2 + i \mu_3) \varepsilon] + \BigO(\varepsilon^2) = 3(\mu_0-\beta_2 \varepsilon) +3 i (\sqrt{7 \mu_0} + \mu_3 \varepsilon)  + \BigO(\varepsilon^2)\\
	{q}_2 &= -10 + 3 \beta_1 \varepsilon + \BigO(\varepsilon^2) \\
	{q}_3 &= 0 \\
	{q}_4 &= -5 \mu_0 -2 (\beta_2+\mu_3 \sqrt{7 \mu_0}) \varepsilon + i [\sqrt{7 \mu_0}(2 + \mu_0) + (2 \mu_3 - \beta_2 \sqrt{7 \mu_0} + \mu_0 \mu_3) \varepsilon ]   +  \BigO(\varepsilon^2) 
\end{aligned}
\]
From \(A q = i \omega_0 q\) we get \(A \overline{q} = -i \omega_0 \overline{q}\).

\vskip 0.2cm
\noindent
Now we must find a vector \(p\) such that \(A^T p = - i \omega_0 \, p\), and also \(A^T \overline{p} = i \omega_0 \, \overline{p}\), that is:

\begin{equation} \label{eq:matriceAtrasposta}
\scalebox{0.8}{$	\begin{pmatrix}
		-2 + i \omega_0 & -\frac{10}{3} + \beta_1 \varepsilon + \BigO(\varepsilon^2) & 0 & \frac{10}{3}  \\ \\
		0 & -\mu_0 + \beta_2 \varepsilon + i \omega_0 + \BigO(\varepsilon^2) & 0 & 50 - \mu_0 + \beta_2 \varepsilon + \BigO(\varepsilon^2) \\ \\
		0 & 0 & \alpha_T \varepsilon  + i \omega_0+ \BigO(\varepsilon^2) & \beta_3 \varepsilon + \BigO(\varepsilon^2) \\ \\
		3 & 0 & 0 & -5 - \beta_3 \varepsilon  + i \omega_0 + \BigO(\varepsilon^2)
	\end{pmatrix}$} \begin{pmatrix}
		\varphi_1 \\ \varphi_2 \\ \varphi_3 \\ \varphi_4
	\end{pmatrix} =  \begin{pmatrix}
		0 \\ 0 \\ 0 \\ 0
	\end{pmatrix}.
\end{equation}

Solving the system, we find
\begin{equation}
	\left\{
	\begin{aligned}
		\varphi_1 & =  \frac {5 + \beta_3 \varepsilon - i \omega_0} {3} \varphi_4 + \BigO(\varepsilon^2), \\
		\varphi_2 & =  \frac {50 - \mu_0+ \beta_2 \varepsilon } {\mu_0- \beta_2 \varepsilon -i \omega_0} \varphi_4 + \BigO(\varepsilon^2)\\
		\varphi_3 &= \frac {- \beta_3 \varepsilon } {\alpha_T \varepsilon  + i \omega_0} \varphi_4 + \BigO(\varepsilon^2)
	\end{aligned}
	\right.
\end{equation}
with \(\varphi_4 \in \mathbb{C}\). Letting
\[
\varphi_4 = 3(-\mu_0 + i \sqrt{7 \mu_0} + (\beta_2 + i \mu_3) \varepsilon) + \BigO(\varepsilon^2)
\]
and recalling
\[
- \frac{3 \beta_3}{\sqrt{7 \mu_0} + 2 \mu_3 \varepsilon} = - \frac{3 \beta_3}{7 \mu_0}(\sqrt{7 \mu_0} - 2 \mu_3 \varepsilon) + \BigO(\varepsilon^2)
\] we obtain

\begin{equation} \label{eq:ptilde}  \begin{aligned}
		{p}_1 = & 2 \mu_0 + i \sqrt{7 \mu_0} (5 + \mu_0) + (\beta_4 + i \beta_5) \varepsilon + \BigO(\varepsilon^2)\\
		{p}_2= & -150 + 3\mu_0 - 3 \beta_2 \varepsilon + \BigO(\varepsilon^2) \\
		{p}_3 = & - \frac {3 \beta_3} {\sqrt{7 \mu_0} + 2 \mu_3 \varepsilon} \left( \sqrt{7 \mu_0}  + i \mu_0 \right) \varepsilon + \BigO(\varepsilon^2) = -  {3 \beta_3}\left( 1 + i \sqrt{\frac {\mu_0} 7}\right) \varepsilon + \BigO(\varepsilon^2)\\
		{p}_4 = & 3(-  \mu_0 + i \sqrt{7 \mu_0}) + 3( \beta_2 + i \mu_3 )  \varepsilon  + \BigO(\varepsilon^2)
	\end{aligned} \end{equation}
where
\[
\begin{aligned}
	\beta_4 &= 2 \mu_3 \sqrt{7 \mu_0} - \beta_3 \mu_0 + 5 \beta_2 \\
	\beta_5 &= 5 \mu_3 + \mu_0 \mu_3 + \beta_3 \sqrt{7 \mu_0} - \beta_2 \sqrt{7 \mu_0}.
\end{aligned}
\]

\noindent
We now proceed by imposing the normalization condition, that is, we look for
 \(\alpha = x + i y \in \mathbb{C}\)
such that
\[ 1= \left\langle \frac p \alpha, q \right\rangle.\]
Recalling that 
\[\left\langle \frac p \alpha, q \right\rangle =\sum_{i=1}^3 \overline{\left( \frac {p_i} \alpha \right)} q_i = \overline {\left(\frac 1 \alpha \right)} \sum_{i=1}^3 \overline{p_i}  q_i = \frac {x+iy}{x^2+y^2} \sum_{i=1}^3 \overline{p_i}  q_i\]

	 we compute
\[\begin{aligned}
	 \overline{p_1} q_1 &=  \left[ 2 \mu_0 - i \sqrt{7 \mu_0} (5 + \mu_0)  + \BigO(\varepsilon)\right]  \left[ 3 \mu_0 +3 i \sqrt{7 \mu_0} + \BigO(\varepsilon)\right] 
	\\ & = 3 \mu_0 (9 \mu_0+35) -3 \,i \mu_0 \sqrt {7 \mu_0}(3+ \mu_0) + \BigO(\varepsilon)
	   \\ \\
	\overline{p_2}  q_2 & = \left[ -150 + 3\mu_0  + \BigO(\varepsilon)\right] \left[ -10 + \BigO(\varepsilon) \right]= {30(50-\mu_0)} + \BigO(\varepsilon)\\ \\
	\overline{p_3}  q_3& =  0\\ \\
	\overline{p_4} q_4& =  \left[-3(  \mu_0 + i \sqrt{7 \mu_0})  + \BigO(\varepsilon) \right] \left[-5 \mu_0  + i \sqrt{7 \mu_0}(2 + \mu_0)   +  \BigO(\varepsilon)  \right]\\ & = 6 \mu_0 (6 \mu_0 +7)  +  3 i \mu_0 \sqrt{7 \mu_0}(3- \mu_0) +\BigO(\varepsilon).
\end{aligned}\]
Summing these contributions, we obtain
\[ \begin{aligned} \left\langle \frac p \alpha, q \right\rangle &= \frac {x+iy} {x^2+y^2} \left[63 \mu_0^2 +117 \mu_0 + 1500 - 6 i \mu_0^2 \sqrt{7 \mu_0} +\BigO(\varepsilon) \right] \\ &= \frac {3} {x^2+y^2} \cdot \\ & \cdot \left[ 21\,x{\mu_0}^{2}+39\,x\mu_0  +500\,x+2\,y \mu_0^2\sqrt {7\mu_0} + i (-2 x \mu_0^2 \sqrt {7 \mu_0}+500y+39y
\mu_0+21y{\mu_0}^{2}) +  \BigO(\varepsilon) 
\right].\end{aligned}\]
Since \(\left\langle \frac{p}{\alpha}, q \right\rangle = 1\), we get the system
\[\left\{
\begin{aligned}
	&  (21 \mu_0^2 + 39 \mu_0 +500)x + 2 \sqrt{7 \mu_0} \mu_0^2 y=  \frac {x^2+y^2} 3, \\
	& -2 \sqrt{7 \mu_0} \mu_0^2 x +(21 \mu_0^2 + 39 \mu_0 +500)y =0.
\end{aligned}
\right.\]
Solving this system, we find
\[y= \frac {2 \sqrt{7 \mu_0} \mu_0^2} {21 \mu_0^2 + 39 \mu_0 +500} x \]
and hence
\begin{equation} \label{eq:sistema} \left\{ \begin{aligned} x & = 3(21 \mu_0^2 + 39 \mu_0+500):= x(\mu_0) \\ \\
	y  & =   {6 \sqrt{7 \mu_0} \mu_0^2}  := y(\mu_0).\end{aligned} \right. \end{equation}
Thus, the normalization factor is given by
\[\alpha(\mu_0)= x(\mu_0) + i y(\mu_0)\]
with \(x(\mu_0), y(\mu_0)\) as in \eqref{eq:sistema}.

\vskip 0.2cm
\noindent
Now we have to compute the inverse matrix \(A^{-1}\).
We have \[
\det A = 10 \alpha_T (\mu_0 - 50) \varepsilon + \BigO(\varepsilon^2)
\] 
and the inverse matrix \(A^{-1}\) is given by

 \[
		A^{-1} =	\frac 1 {\det A} \scalebox{0.8}{$ \begin{pmatrix} 
			5 \alpha_T \mu_0 \varepsilon+ \BigO(\varepsilon^2) & 3 \alpha_T(50- \mu_0) \varepsilon + \BigO(\varepsilon^2) & - 3 \beta_3 \mu_0 \varepsilon + \BigO(\varepsilon^2)&  3 \alpha_T \mu_0 \varepsilon + \BigO(\varepsilon^2)\\ \\
			- \frac {50} 3 \alpha_T \varepsilon + \BigO(\varepsilon^2) &  \BigO(\varepsilon^2) &  10 \beta_3 \varepsilon + \BigO(\varepsilon^2) & -10 \alpha_T \varepsilon + \BigO(\varepsilon^2)\\ \\
			0 & 0 & *_{3,3}  &0 \\ \\
			\frac {20} 3 (\mu_0 -25) \alpha_T \varepsilon+ \BigO(\varepsilon^2) &2 \alpha_T (50 - \mu_0 ) \varepsilon + \BigO(\varepsilon^2) & -2\beta_3 \mu_0 \varepsilon + \BigO(\varepsilon^2) & 2 \alpha_T \mu_0 \varepsilon + \BigO(\varepsilon^2)
		\end{pmatrix} $} \]
	with
	\[
	*_{3,3} = 10 \mu_0 - 500 + (-3 \beta_1 \mu_0 - 2 \beta_3 \mu_0 + 150 \beta_1 - 10 \beta_2) \varepsilon + \BigO(\varepsilon^2)
	\]
	and then 
	\begin{equation} \label{eq:matriceA_inversa}
		A^{-1} = \scalebox{0.8}{$ \begin{pmatrix} 
			\displaystyle \frac {\mu_0} {2(\mu_0-50)}+ \BigO(\varepsilon) &\displaystyle  - \frac 3 {10} + \BigO(\varepsilon) & \displaystyle - \frac {3 \beta_3 \mu_0} {10 \alpha_T (\mu_0 -50)}  + \BigO(\varepsilon)&  \displaystyle \frac {3 \mu_0} {10(\mu_0 -50)}  + \BigO(\varepsilon)\\ \\
			\displaystyle	- \frac {5} {3(\mu_0-50)} + \BigO(\varepsilon) &  \BigO(\varepsilon) &  \displaystyle \frac {\beta_3} {\alpha_T (\mu_0-50)} + \BigO(\varepsilon) & \displaystyle - \frac 1 {\mu_0-50} + \BigO(\varepsilon)\\ \\
			0 & 0 &  \frac1 {\alpha_T  \varepsilon } - \frac {3 \beta_1(\mu_0 -50) +10 \beta_2 +  2 \beta_3 \mu_0} {10 \alpha_T (\mu_0-50)} + \BigO(\varepsilon) &0 \\ \\
			\displaystyle \frac {2(\mu_0-25)} {3 (\mu_0-50)}+ \BigO(\varepsilon) &\displaystyle - \frac 15 + \BigO(\varepsilon) &\displaystyle -\frac {\beta_3 \mu_0} {5 \alpha_T(\mu_0 -50)} + \BigO(\varepsilon) & \displaystyle \frac {\mu_0} {5(\mu_0-50)} + \BigO(\varepsilon)
		\end{pmatrix}. $}
\end{equation}

We note that now all entries of the matrix \(A^{-1}\) are of the form \(const + \BigO(\varepsilon)\), except for the entry \(*_{3,3} = \BigO(\varepsilon^{-1})\).

Following Kuznetsov \cite{KuznetsovFourth}, we now compute the multilinear functions \(B(x, y)\) and \(C(x, y, z)\), given in \eqref{eq:B} and \eqref{eq:C} 
where \(F(\xi_1,\xi_2 ,\xi_3)= (F_1, F_2, F_3)(\xi_1,\xi_2,\xi_3)\)
denotes the vector field 
\[\begin{aligned}
	F_1(\xi_1, \xi_2, \xi_3, \xi_4) & = - 2 \xi_1 +3 \xi_4 \\
	F_2(\xi_1, \xi_2, \xi_3, \xi_4) & = 50 \xi_2 - \frac 12 \xi_2^2 - \delta_0 \xi_1 \xi_2 \\
	F_3(\xi_1, \xi_2, \xi_3, \xi_4) & = - \varepsilon  \xi_1 \xi_3  + k \varepsilon \frac{\xi_1 \xi_4}{\varepsilon \xi_1 \xi_4 + 1} \xi_3\\
	F_4(\xi_1, \xi_2, \xi_3, \xi_4) & = - 5 \xi_4 + \delta \xi_1 \xi_2 + \varepsilon \xi_1 \xi_3 - \varepsilon \xi_1 \xi_4,\\
\end{aligned}\]
with \(\delta_0 = \displaystyle \frac 5 {3 \mu_0}\). It is immediate that all third-order derivatives vanish, except for those involving \(F_3\).  
In particular, for \(F_1, F_2,\) and \(F_4\) we have 
\( C_i(x, y, z) \equiv 0 \) for \(i=1,2,4\).  

Moreover,
\[
\frac{\partial^3 F_3}{\partial \xi_1 \partial \xi_3 \partial \xi_4} 
= \frac{k \varepsilon}{1 + \varepsilon x w} + \BigO(\varepsilon^2) 
= k \varepsilon + \BigO(\varepsilon^2),
\]
while all the other third-order derivatives of \(F_3\) are of order \(\BigO(\varepsilon^2)\).  
Therefore,
\[
C_3(x,y,z) 
= k \varepsilon \,(x_1 y_3 z_4 + x_1 y_4 z_3 + x_3 y_1 z_4 
+ x_3 y_4 z_1 + x_4 y_1 z_3 + x_4 y_3 z_1) + \BigO(\varepsilon^2) 
= \BigO(\varepsilon).
\]

\noindent
To compute \(B\), we need to evaluate the second-order derivatives. All second derivatives of \(F_1\) vanish, so \( B_1(x, y) \equiv 0 \).  
To compute \( B_2, \, B_3 \) and \( B_4 \), recall the definition
\[B_i(x,y)= \sum_{j,k=1}^3 \left( \frac {\partial^2 F_i} {\partial \xi_j \partial \xi_k} \right) x_j y_k=  \left \langle D^2 F_i x , y \right \rangle\]
where \(D^2F_i\) is the Hessian matrix of \(F_i, i=2,3,4\).
For \(i=2\) we have
\[ \frac {\partial F_2} {\partial \xi_1} = - \delta_0 \xi_2 , \quad \frac {\partial F_2} {\partial \xi_2} = 50 - \delta_0 \xi_1 - \xi_2 , \quad \frac {\partial F_2} {\partial \xi_3} =  \frac {\partial F_2} {\partial \xi_4} = 0. \]
Thus, the Hessian matrix of \(F_2\) is:
\[D^2 F_2 = D^2 F_2(P_0) = \begin{pmatrix} 
	0 & -\delta_0 & 0 & 0\\ 
	-\delta_0 & -1 & 0 & 0\\
	0 & 0 & 0 & 0 \\ 
	0& 0 & 0 & 0  
\end{pmatrix} \]
so
\[B_2(x,y) = -  \delta_0 x_1y_2 - \delta_0 x_2 y_1 -x_2 y_2.\]
For \(i=3\), we have:
\[
\begin{aligned}
	\frac{\partial F_3}{\partial \xi_1} &= -\varepsilon \xi_3 
	+ \frac{k \varepsilon \xi_4}{(\varepsilon \xi_1 \xi_4 + 1)^2}\,\xi_3, 
	&\qquad 
	\frac{\partial F_3}{\partial \xi_2} &= 0, \\[6pt]
	\frac{\partial F_3}{\partial \xi_3} &= -\varepsilon \xi_1 
	+ k \varepsilon \frac{\xi_1 \xi_4}{\varepsilon \xi_1 \xi_4 + 1}, 
	&\qquad 
	\frac{\partial F_3}{\partial \xi_4} &= \frac{k \varepsilon \xi_1}{(\varepsilon \xi_1 \xi_4 + 1)^2}\,\xi_3 .
\end{aligned}
\]
The Hessian matrix is given by
\begin{equation} \label{eq:hessianaf3}
	D^2 F_3 =
	\begin{pmatrix} 
		-2 k \varepsilon^2 \dfrac{\xi_3 \xi_4^2}{(\varepsilon \xi_1 \xi_4 + 1)^3} 
		& 0 
		& - \varepsilon + \dfrac{k \varepsilon \xi_4}{(\varepsilon \xi_1 \xi_4 +1)^2} 
		& k \varepsilon \xi_3 \dfrac{1 - \varepsilon \xi_1 \xi_4}{(\varepsilon \xi_1 \xi_4 + 1)^3} \\[8pt]
		0 & 0 & 0 & 0 \\[8pt]
		- \varepsilon + \dfrac{k \varepsilon \xi_4}{(\varepsilon \xi_1 \xi_4 +1)^2} 
		& 0 
		& 0 
		& \dfrac{k \varepsilon \xi_1}{(\varepsilon \xi_1 \xi_4 +1)^2} \\[8pt]
		k \varepsilon \xi_3 \dfrac{1 - \varepsilon \xi_1 \xi_4}{(\varepsilon \xi_1 \xi_4 + 1)^3} 
		& 0 
		& \dfrac{k \varepsilon \xi_1}{(\varepsilon \xi_1 \xi_4 +1)^2} 
		& -2 k \varepsilon^2 \dfrac{\xi_1^2 \xi_3}{(\varepsilon \xi_1 \xi_4 + 1)^3}
	\end{pmatrix}.
\end{equation}

Evaluating this matrix at the point 
\[
\xi_1 = x_0, \quad 
\xi_2 = y_0, \quad 
\xi_3 = z_0 = 0, \quad 
\xi_4 = w_0 = \tfrac{2}{3} x_0,
\]
and using the expansion
\[
\frac{1}{(\varepsilon x w + 1)^2} 
= 1 - 2 x w \varepsilon + \BigO(\varepsilon^2),
\]
we find
\[
\begin{aligned}
	- \varepsilon + \frac{k \varepsilon w}{(\varepsilon x w +1)^2} 
	&= (kw -1) \varepsilon  + \BigO(\varepsilon^2) \\
	&= \left(\tfrac{2}{3} k x_0 -1 \right) \varepsilon + \BigO(\varepsilon^2) \\
	&= \left(\tfrac{2}{3} k x_T -1 \right) \varepsilon + \BigO(\varepsilon^2).
\end{aligned}
\]

Thus, the Hessian becomes
\[
D^2 F_3(P_0) =
\begin{pmatrix} 
	0 & 0 & \left(\tfrac{2}{3} k x_T -1 \right) \varepsilon + \BigO(\varepsilon^2) & 0 \\[8pt]
	0 & 0 & 0 & 0 \\[8pt]
	\left(\tfrac{2}{3} k x_T -1 \right) \varepsilon + \BigO(\varepsilon^2) & 0 & 0 & k x_T \varepsilon + \BigO(\varepsilon^2) \\[8pt]
	0 & 0 & k x_T \varepsilon + \BigO(\varepsilon^2) & 0
\end{pmatrix}.
\]

Consequently,
\[
\begin{aligned}
	B_3(x,y) &= \Big[ 
	\left(\tfrac{2}{3} k x_T -1 \right) (x_1 y_3 + x_3 y_1) 
	+ k x_T (x_3 y_4 + x_4 y_3)
	\Big] \varepsilon + \BigO(\varepsilon^2) \\
	&= \BigO(\varepsilon).
\end{aligned}
\]

Finally, for \(i=4\) we have
\[
\frac{\partial F_4}{\partial \xi_1} = \delta_0 \xi_2 + \varepsilon \xi_3 - \varepsilon \xi_4, 
\quad 
\frac{\partial F_4}{\partial \xi_2} = \delta_0 \xi_1, 
\quad 
\frac{\partial F_4}{\partial \xi_3} = \varepsilon \xi_1, 
\quad 
\frac{\partial F_4}{\partial \xi_4} = -5 - \varepsilon \xi_1.
\]

The corresponding Hessian matrix is
\[
D^2 F_4 = D^2 F_4(P_0) =
\begin{pmatrix}
	0 & \delta_0 & \varepsilon & -\varepsilon \\
	\delta_0 & 0 & 0 & 0 \\
	\varepsilon & 0 & 0 & 0 \\
	-\varepsilon & 0 & 0 & 0
\end{pmatrix}.
\]

Therefore,
\[
\begin{aligned}
	B_4(x,y) &= \delta_0 x_1 y_2 + \varepsilon x_1 y_3 - \varepsilon x_1 y_4 
	+ \delta_0 x_2 y_1 + \varepsilon x_3 y_1 - \varepsilon x_4 y_1 \\
	&= \delta_0 (x_1 y_2 + x_2 y_1) + \BigO(\varepsilon).
\end{aligned}
\]

To compute the coefficient \(	a(\mu(\varepsilon),\varepsilon ) \) following Kuznetsov's formulas, we begin with the evaluation of \(B(q, \overline{q})\). We now consider a first-order simplified expression for \(q\) and for its complex conjugate \(\overline{q}\), namely
\[\begin{aligned} q & =\left( 3 (\mu_0 + i \sqrt{7 \mu_0}) + \BigO(\varepsilon), \;  -10 + \BigO(\varepsilon),\;  0, \;-5 \mu_0 + i \sqrt{7 \mu_0}(2 + \mu_0)+ \BigO(\varepsilon) \right)  \\
\overline{q} & = \left( 3 (\mu_0 - i \sqrt{7 \mu_0})+ \BigO(\varepsilon), \;  -10+ \BigO(\varepsilon),\;  0, \;-5 \mu_0 - i \sqrt{7 \mu_0}(2 + \mu_0)+ \BigO(\varepsilon) \right).\end{aligned}\] 
The components of the vector \(B(q,\overline{q})\) are given by
\[
\begin{aligned}
	B_1(q,\overline{q}) &= 0, \\
	 B_2(q,\overline{q}) & = -\delta_0 q_1 \overline{q}_2 - \delta_0 q_2 \overline{q}_1 - q_2 \overline{q}_2 
	=  \BigO(\varepsilon), \\ 
	B_3(q,\overline{q}) &= \Big[ 
	\left(\tfrac{2}{3} k x_T -1 \right) (q_1 \overline{q}_3 + q_3 \overline{q}_1) 
	+ k x_T (q_3 \overline{q}_4 + q_4 \overline{q}_3)
	\Big] \varepsilon + \BigO(\varepsilon^2) =  0,  \\
	 B_4(q,\overline{q}) &= \delta_0 q_1 \overline{q}_2 + \delta_0 q_2 \overline{q}_1 + \BigO(\varepsilon) 
	= -100 + \BigO(\varepsilon).
\end{aligned}
\]

\noindent
Now we compute \( v= A^{-1} B(q, \overline{q})\) with \(A^{-1}\) as in \eqref{eq:matriceA_inversa} and with straightforward computations we obtain
\[ v= A^{-1} B(q, \overline{q})  = \left( - \frac {30 \mu_0}{\mu_0-50} + \BigO(\varepsilon), \; \frac {100} {\mu_0 -50}+ \BigO(\varepsilon), \; 0,  \; - \frac {20 \mu_0} {\mu_0 -50}+ \BigO(\varepsilon)\right).\]
We note that the entry in position \((3,3)\) of the matrix \(A^{-1}\) is \(\BigO(\varepsilon^{-1})\), and therefore it is  important that  the coefficient of \(\varepsilon\) in \(B_3(q, \overline{q})\) is indeed zero.

We now compute the vector
\(B\bigl(q, A^{-1} B(q, \overline{q})\bigr) = B(q,v)\), 
and obtain
\[
\begin{aligned}
	B_1(q,v) &= 0, 
	&\qquad B_2(q,v) &= -\delta_0 q_1 v_2 - \delta_0 q_2 v_1 - q_2 v_2 
	= - i \, \frac{500 \sqrt{7 \mu_0}}{\mu_0(\mu_0-50)} + \BigO(\varepsilon), \\[6pt]
	B_3(q,v) &= 0, 
	&\qquad B_4(q,v) &= \delta_0 q_1 v_2 + \delta_0 q_2 v_1 + \BigO(\varepsilon) 
	= \frac{500(2 \mu_0 + i \sqrt{7 \mu_0})}{\mu_0(\mu_0-50)} + \BigO(\varepsilon).
\end{aligned}
\]

\noindent
The inner sum evaluates to
\[ \begin{aligned}
\sum_{k=1}^4 \overline{p_k} B_k(q,v)  & = \overline{p_2} B_2(q,v) 
 + \overline{p_4} B_4(q,v) \\
& = -1500 \, \frac{2\mu_0^2 - 7\mu_0 + i \sqrt{7\mu_0}(4\mu_0 -50)}
{\mu_0(\mu_0-50)} + \BigO(\varepsilon),
\end{aligned}\]
hence,
\[ \displaystyle
\left \langle \frac{p}{\alpha}, B(q, v) \right \rangle
= - \frac{1500(x+iy)}{x^2+y^2}
\left[
\frac{2\mu_0^2 - 7\mu_0 + i \sqrt{7\mu_0}(4\mu_0 -50)}
{\mu_0(\mu_0-50)}
\right] + \BigO(\varepsilon).
\]

So the real part is:
\[ \begin{aligned} &\operatorname{Re} \left( \left \langle \frac p \alpha, B(q, A^{-1}B(q, \overline{q})) \right \rangle\right) = \operatorname{Re} \left( \left \langle \frac p \alpha, B(q,v) \right \rangle \right) \\
	& = - \frac {1500} {\mu_0 (\mu_0-50) (x^2+y^2)} \operatorname{Re} \left( (x+iy)\left[  {2 \mu_0^2 - 7 \mu_0 + i \sqrt{7 \mu_0}(4 \mu_0 -50)}\right]   \right) + \BigO(\varepsilon) \\
	& = - \frac {1500 \, [(2 \mu_0^2- 7 \mu_0) x  - \sqrt{7 \mu_0} y (4 \mu_0 -50)]} {\mu_0(\mu_0-50) (x^2+y^2)}  + \BigO(\varepsilon).\end{aligned}\]
To compute the second term of \( 	a(\mu(\varepsilon),\varepsilon )  \), we evaluate \( B(q, q) \). We have 
\[B(q,q)=\left (0, \;  i 100 \sqrt{\frac 7 {\mu_0}} + \BigO(\varepsilon), \; 0 , \; -100  -i \, 100  \sqrt{\frac 7 {\mu_0}} + \BigO(\varepsilon) \right ).\]

\noindent
We now compute the vector
\(w=(2 i \omega_0 I -A)^{-1}B(q,q)\).
We have
\[ -A + i 2 \omega_0 I =
\begin{pmatrix}
	2 + i2 \sqrt{7 \mu_0} + \BigO(\varepsilon) & 0 & 0 & -3 \\ \\
	\frac{10}{3}+ \BigO(\varepsilon) & \mu_0 + i2 \sqrt{7 \mu_0} + \BigO(\varepsilon) & 0 & 0 \\ \\
	0 & 0 & i2 \sqrt{7 \mu_0} + \BigO(\varepsilon) & 0 \\ \\
	-\frac{10}{3} & -50 + \mu_0 + \BigO(\varepsilon) & + \BigO(\varepsilon) & 5 + i2 \sqrt{7 \mu_0} + \BigO(\varepsilon)
\end{pmatrix}.
\]
The determinant is:
\begin{equation} \label{eq:deTinvAi}
	\det(-A + 2i \omega_0 I) = 588 \mu_0^2 + 4i \sqrt{7 \mu_0} (-14 \mu_0^2 - 103 \mu_0 + 250) + \BigO(\varepsilon) = D + \BigO(\varepsilon),
\end{equation}
with
\begin{equation} \label{eq:D}
	D:= 588 \mu_0^2 + 4i \sqrt{7 \mu_0} (-14 \mu_0^2 - 103 \mu_0 + 250) .
\end{equation}
The inverse matrix is then:
\begin{equation} \label{eq:invAi}
	(-A + 2i \omega_0 I)^{-1} = \frac{1}{\det(-A + 2i \omega_0 I)} T
\end{equation}
where  \(T\) is given by:

\[
\resizebox{\textwidth}{!}{$
	T =
	\begin{pmatrix}
		- 28 \mu_0 (5 + \mu_0) - 46 i \mu_0 \sqrt{7 \mu_0} + \BigO(\varepsilon) &
		6 i \sqrt{7 \mu_0}(-\mu_0 + 50) + \BigO(\varepsilon) &
	 \BigO(\varepsilon) &
		- 6 \mu_0 (14 - i \sqrt{7 \mu_0}) + \BigO(\varepsilon) \\ \\
		\frac{280}{3} \mu_0 - \frac{100}{3} i \sqrt{7 \mu_0} + \BigO(\varepsilon) &
		- 196 \mu_0 - 56 i \mu_0 \sqrt{7 \mu_0} + \BigO(\varepsilon) &
		 \BigO(\varepsilon) &
		- 20 i \sqrt{7 \mu_0} + \BigO(\varepsilon) \\ \\
		\BigO(\varepsilon^2) &
		\BigO(\varepsilon^2) &
		- 28 \mu_0^2 - 206 \mu_0 + 500 - 42 i \mu_0 \sqrt{7 \mu_0} + \BigO(\varepsilon) &
		\BigO(\varepsilon^2) \\ \\
		- \frac{280}{3} \mu_0 - \frac{40}{3} i \sqrt{7 \mu_0} (25 - \mu_0) + \BigO(\varepsilon) &
		28 \mu_0^2 - 1400 \mu_0 - 4 i \sqrt{7 \mu_0} (\mu_0 - 50) + \BigO(\varepsilon) &
		 \BigO(\varepsilon) &
		- 28 \mu_0 (\mu_0 + 2) - 52 i \mu_0 \sqrt{7 \mu_0} + \BigO(\varepsilon)
	\end{pmatrix}.
	$}
\]
For the reader's convenience, we rewrite below the matrix \(T\) in a more readable form, first showing the first two columns and then the third and fourth ones:
\[
T =
\left(
\begin{array}{cc}
	- 28 \mu_0 (5 + \mu_0) - 46 i \mu_0 \sqrt{7 \mu_0} + \BigO(\varepsilon) &
	6 i \sqrt{7 \mu_0}(-\mu_0 + 50) + \BigO(\varepsilon) \\[6pt]
	\frac{280}{3} \mu_0 - \frac{100}{3} i \sqrt{7 \mu_0} + \BigO(\varepsilon) &
	- 196 \mu_0 - 56 i \mu_0 \sqrt{7 \mu_0} + \BigO(\varepsilon) \\[6pt]
	\BigO(\varepsilon^2) &
	\BigO(\varepsilon^2) \\[6pt]
	- \frac{280}{3} \mu_0 - \frac{40}{3} i \sqrt{7 \mu_0} (25 - \mu_0) + \BigO(\varepsilon) &
	28 \mu_0^2 - 1400 \mu_0 - 4 i \sqrt{7 \mu_0} (\mu_0 - 50) + \BigO(\varepsilon)
\end{array}
\right. \]
\vskip 0.2cm
\[\left.
\begin{array}{cc}
	 \BigO(\varepsilon) &
	- 6 \mu_0 (14 - i \sqrt{7 \mu_0}) + \BigO(\varepsilon) \\[6pt]
	 \BigO(\varepsilon) &
	- 20 i \sqrt{7 \mu_0} + \BigO(\varepsilon) \\[6pt]
	- 28 \mu_0^2 - 206 \mu_0 + 500 - 42 i \mu_0 \sqrt{7 \mu_0} + \BigO(\varepsilon) &
	\BigO(\varepsilon^2) \\[6pt]
	 \BigO(\varepsilon) &
	- 28 \mu_0 (\mu_0 + 2) - 52 i \mu_0 \sqrt{7 \mu_0} + \BigO(\varepsilon)
\end{array}
\right).
\]

\noindent
The components of \(w=(2 i \omega_0 I -A)^{-1}B(q,q) \) are the following
\[\begin{aligned}
	w_1 &= \frac {600} {D} \left[14(-25+2 \mu_0) +i\sqrt {7 \mu_0}(14-\mu_0)\right] + \BigO(\varepsilon),\\
	w_2 & = \frac {400} {D}(-35+98\,\mu_0 -44\,i\sqrt {7 \mu_0}) + \BigO(\varepsilon),\\
	w_3 &=  0, \\
	w_4 &= \frac {400} {D}	\left[7(-50-10\mu_0 + \mu_0^2) +3\,i\sqrt {7 \mu_0}(9 \mu_0-112)
	\right] + \BigO(\varepsilon).
\end{aligned}\]

\noindent
We now compute the vector \( B(\overline{q}, w) \):
\[\begin{aligned}
	B_1(\overline{q}, w) &= 0, \\
	B_2(\overline{q},w) & = \frac {14000} {\mu_0 D} \left[{14 \mu_0^2 + 59 \mu_0 -250 +i\sqrt{7 \mu_0}(7\mu_0+5)} \right] +\BigO(\varepsilon),\\
	B_3(\overline{q},w) & = 0, \\
	B_4(\overline{q},w) & = - \frac {2000} {\mu_0 D} \left[{-98 \mu_0^2+ 483 \mu_0 -1750 +i\sqrt{7 \mu_0}(137\mu_0+35)} \right] +\BigO(\varepsilon).
\end{aligned}\]

\noindent
Then we recall that
\[ \left \langle \frac p \alpha, B(\overline{q}, (2i \omega_0 I-A)^{-1}B(q,q )) \right \rangle  =  \left \langle \frac p \alpha, B(\overline{q},w) \right \rangle   = \frac {x + i y} {x^2+y^2} \; \sum_{k=1}^4  \overline{p_k} B_k(\overline{q},w)\]
and we obtain
\[\begin{aligned}  \sum_{k=1}^4  \overline{p_k} B_k(\overline{q},w) 
	 & =\frac {6000} {\mu_0 D} \left[7(-709\mu_0^2 - 3485\mu_0+12500)  + i \sqrt{7 \mu_0} (88 \mu_0^2 -1897 \mu_0 -3500)\right] +\BigO(\varepsilon) \\ 
	& = \frac {6000} {\mu_0 D} \left[E+ i \sqrt{7 \mu_0}F\right] +\BigO(\varepsilon)
\end{aligned}\]
with \(D\) as in \eqref{eq:D} and
\[
E  := 7(-709\mu_0^2 - 3485\mu_0+12500), \qquad  
	F :=88 \mu_0^2 -1897 \mu_0 -3500.
\]
Therefore,
\[\left \langle \frac p \alpha, B(\overline{q}, (2i \omega_0 I-A)^{-1}B(q,q )) \right \rangle  =   \frac {6000(x + i y)} {(x^2+y^2)\mu_0 D}  \left[E + i \sqrt{7 \mu_0} F\right] +\BigO(\varepsilon).\]	

\noindent
We find
\[\begin{aligned}  \operatorname{Re} \left(
	\left \langle \frac p \alpha, B(\overline{q}, (2i \omega_0 I-A)^{-1}B(q,q )) \right \rangle \right)& =   \frac {6000} {(x^2+y^2)\mu_0 } \operatorname{Re} \left(\frac {(x+iy) (E+ i \sqrt{7 \mu_0} F) } {D}\right) +\BigO(\varepsilon)\\
	& = \frac {6000} {(x^2+y^2)\mu_0 } \operatorname{Re} \left( \frac { x E- y\sqrt{7 \mu_0} F  + i (x \sqrt{7 \mu_0} F+ yE )} { D}\right) +\BigO(\varepsilon) \\
	&= \frac {6000} {(x^2+y^2)\mu_0 }  \operatorname{Re} \left( z \right) +\BigO(\varepsilon)
\end{aligned}\]	
where
\(z = \dfrac{\ell_1 + i k_1}{\ell_2 + i k_2}\) with
\[
\begin{aligned}
	\ell_1 &:= x E - y\sqrt{7 \mu_0}\, F, \qquad \qquad&
	k_1 &:= x \sqrt{7 \mu_0}\, F + yE,\\
	\ell_2 &:= 588 \mu_0^2, &
	k_2 &:= 4 \sqrt{7 \mu_0}\, (-14 \mu_0^2 - 103 \mu_0 + 250).
\end{aligned}
\]
From the identity
\(\operatorname{Re}(z) = \frac{\ell_1 \ell_2 + k_1 k_2}{\ell_2^2 + k_2^2}\),
it follows that
\[\begin{aligned}
	\operatorname{Re}(z)
	& = \sqrt {7 \mu_0} \frac { 8078\,{\mu_0}^{4}+161654\,{\mu_0}^{3}+80205\,{\mu_0}^{2}-2158750 \mu_0+3125000}
	{ 4\mu_0\, \left( 196
		\,{\mu_0}^{4}+5971\,{\mu_0}^{3}+3609\,{\mu_0}^{2}-51500\,\mu_0+62500
		\right) } y \\ & - \frac {1232\,{\mu_0}^{4}+86729\,{\mu_0}^{3}+245904\,{\mu_0}^{2}-1723750\,{\mu_0}+875000}{ 4\, \left( 196
		\,{\mu_0}^{4}+5971\,{\mu_0}^{3}+3609\,{\mu_0}^{2}-51500\,\mu_0+62500
		\right) }x.
\end{aligned}\]
\noindent
Finally, recalling \( \omega_0 = \sqrt{7 \mu_0} + \BigO(\varepsilon)\), we find

\[\begin{aligned} a(\mu(\varepsilon),\varepsilon) & = \frac 1 {2 \omega_0} \operatorname{Re}
	\left( \left \langle \frac p \alpha, C(q,q, \overline{q}) \right \rangle - 2 \left \langle \frac p \alpha, B(q, A^{-1}B(q, \overline{q}))\right \rangle + \left \langle \frac p \alpha, B(\overline{q}, (2 i \omega_0 I - A)^{-1} B(q,q))\right \rangle\right) \\
	& = - \frac 1 {\omega_0} \operatorname{Re}  \left( \left \langle \frac p \alpha, B(q, A^{-1}B(q, \overline{q}))\right \rangle \right) + \frac 1 {2 \omega_0} \operatorname{Re} \left( \left \langle \frac p \alpha, B(\overline{q}, (2 i \omega_0 I - A)^{-1} B(q,q))\right \rangle \right) + \BigO(\varepsilon)\\
	&=  \frac 1 {\omega_0} \left\{ \frac {1500} {x^2+y^2} \left[ \frac {(2 \mu_0^2- 7 \mu_0) x  - \sqrt{7 \mu_0} y (4 \mu_0 -50)}  {\mu_0(\mu_0-50)}\right] +  \frac {3000 \operatorname{Re} (z)} {(x^2+y^2)\mu_0 } \right\} + \BigO(\varepsilon).
\end{aligned}\]
Now recall that \( \mu_0 = \displaystyle \frac{-59 + \sqrt{17481}}{14} , \; x = 3(21 \mu_0^2 + 39 \mu_0+500), \; 	y  =   {6 \sqrt{7 \mu_0} \mu_0^2}\), hence,
substituting into the expression gives 
\[a(\mu(\varepsilon),\varepsilon) = a_0  + \BigO(\varepsilon),\]
with \(a_0 \sim -0.04869322966\).

\bibliographystyle{plain} 
\bibliography{references} 
\end{document}